\documentstyle[12pt]{article}

\newcommand{\argm}{\mathop{\rm argmin}\limits}

\begin{document}
%\begin{frontmatter}
%\title{Universal Adaptive Estimations and Confidence Intervals in the Nonparametric 
%Statistics.}
\begin{center}
{\bf UNIVERSAL ADAPTIVE ESTIMATIONS AND CONFIDENCE INTERVALS IN THE NONPARAMETRIC STATISTICS.}\\
\vspace{3mm}

 $ {\bf E.Ostrovsky^a, L.Sirota^b } $ \\
\vspace{2mm}
$ ^a $ Corresponding \ Author. \  {\it Department of Mathematics, Ben Gurion University of Negev,
84105, Beer - Sheva, Ben Gurion street, 2, ISRAEL. }\\
\end{center}
E - mail: \ galaostr@cs.bgu.ac.il\\
\begin{center}
$ ^b $ {\it  Department \ of \ Mathematics. 59200, Bar - Ilan University, Ramat - Gan, ISRAEL.} \\
\end{center}
E - mail:  \ sirota@zahav.net.il \\

%\author[Ostrovsky]{E.Ostrovsky},
%\author[Sirota]{L.Sirota}
%\address[Ostrovsky]{ Corresponding Author. Department of Mathematics,
%Ben Gurion University of Negev, 84105, Beer - Sheva, Ben Gurion street, 2, Israel; \\
% E - mail: \ galaostr@cs.bgu.ac.il}
%\address[Sirota]{Department of Mathematics, 52900, Bar - Ilan University, Ramat - Gan, 
%Israel; \\
%E - mail: \ sirota@zahav.net.il }

\normalsize

\vspace{3mm}

%\begin{abstract}

 {\bf Abstract}. The paper considers so-called adaptive estimations of regression, 
distribution density and spectral density of a Gaussian stationary sequence, 
asymptotically optimal in order at a growing number of observation on any 
regular subspace compactly embedded in space $ L_2 $, and confidence intervals, 
also adaptive, are constructed on their basis for the estimated functions in 
an integral norm.

%\end{abstract}

%\begin{keyword}
\vspace{2mm}


\section{Introduction.} \par

  Among the latest fashions in nonparametric statistics are the so-called 
adaptive estimations (AE), i.e. estimations that use no a priori information 
about the estimated function. Many publications have recently appeared where 
AE are constructed which are optimal in order at a growing number of current 
observations on a continuum of various functional classes (cf. References for 
a list of works on AE, which does not, however, claim to be exhaustive).\par
 In (Polyak B., at al., 1990), (Polyak B. at al., 1992), (Golubev G. at al., 1992)
for instance, AE were constructed for the problem of 
estimating regression (R) which are optimal in order on many subspaces of 
space $ L_2 $, and non-adaptive confidence intervals were elaborated on the 
basis of the obtained estimations for the estimated regression function also 
in norm $ L_2 $, which later were somewhat improved in (Golubev at al., 1992).\par

 In  (Efroimovich S., 1985)  AE were 
constructed for problem $(D)$ of estimating distribution density, which are 
optimal on ellipsoids in $L_2 $, while in (Golubev, 1994) AE were constructed for 
problem $ (S) $ of estimating spectral (smooth) density, and so on.\par
 In numerous publications by D. Donoho et al. (Donoho D at al., 1993(1), 1993(2),
1996,  1999(1), 1999(2) )  and in some others AE are 
constructed (and implemented) which are optimal in order on a number of Besov 
spaces. In those papers as well as in  (Golybev G. at al., 1994), (Nussbaum M., 1985),
(Tony Cai at al., 1999), (Lee G., 2003)  diverse orthonormalized 
systems of functions are used to construct AE, 
such as wavelets, wedgelets, unconditional bases, splines, Demmler - Reinsch 
bases, Ridgelets  (Candes E.J., 2003), (Dette H., 2003) etc. \par

 The recent results about kernel estimations in the considered problems see,
for example, (AAD W Van Der Vaart at al., 2003), Allal J., at al.,
2003), (Corinne Berzin at al., 2003).\par

  In our present work, as in the previous ones  (Ostrovsky E.I., 1996, 1997(1); 1997(2),
1999)  we construct and analyze AE on the basis of the classical apparatus of the 
well known trigonometric  approximation theory (Nikolsky S., 1951), (Timan A., 1960),
(Bernstein S, 1952, 1954). \par 
 
{\bf The AE proposed herein feature a speed of convergence  which is optimal in 
order on any regular subspace compactly embedded in space $L_2 $, the 
estimations are universal and very simple in form, which significantly 
facilitates their implementation; finally, we construct exponential adaptive 
confidence intervals (ACI), i.e. such that the tail of the confidence 
probability decreases with exponential speed.}\par
 To the best of our knowledge, adaptive confidence intervals first appeared 
in our publications (Ostrovsky E. at al., 1996, 1997, 1999). 
The precursor for the present paper is perhaps (Ostrovsky E. at al., 1997); in 
comparison with it we now improve the confidence interval 
and strengthen the convergence type of random values - instead of convergence 
by probability we establish convergence with unit probability; {\it (we 
stipulate here that all convergences of a random values sequence are 
understood with probability 1 only.) } \par

\section{Problem statement. Denotations. Conditions.}

 The following three problems are  classical in  nonparametric statistics.\par
{\bf R.\ The regression problem.} Let $ f(x), \ x \in [0,1] $ be an unknown function, 
Riemann-integrable with a square and measured at points of net 
$ x_i = x_{i,n}= i/n, \ i = 1,2,\ldots,n $ with random independent centered 
identically distributed 
errors $ \{\xi_i\}: \ y_i = f(x_i) + \xi_i. $ 
It is required to estimate the function $ f(x) $ with the best possible 
precision from the values $ \{y_i\} $.\par
{\bf D.\ Estimation of distribution density.} On the basis of a set of 
independent identically distributed values $\{\xi_i\}, \ \xi_i \in [0,1], \ 
i = 1,2,\ldots,n $ it is required to estimate their common density $f(x) $ 
(assumed to exist).\par
{\bf S. \ Spectral statistics.} Let $\{\xi_i\} $ be a Gaussian stationary centered 
sequence with spectral density $f(x). $ The estimation object is $f(x).$ 
 We assume for convenience that $ x \in [0,1]. $\par
 It is supposed that all the estimated functions $ f(\cdot) \in L_2[0,1], $ 
therefore they are expanded in the norm of this space into a Fourier series in 
the complete orthonormalized trigonometric system $\{\varphi_j(\cdot)\} $ on 
set [0, 1]: $ \varphi_1(x) = 1; $
$$ 
l > 1 \Rightarrow  \varphi_{2l}(x) = \sqrt{2} \cos(2\pi l x); \ \ 
\varphi_{2l+1}(x) = \sqrt{2} \sin(2\pi l x);
$$
$$
f(x) = \sum_{j=1}^{\infty} c_j \varphi_j(x); \ \ c_j = \int_0^1 \varphi_j(x)
f(x) dx.
$$
 Let us set $ \rho(N) = \rho(f,N) = \sum_{j=N+1}^{\infty} c_j^2. $ 
Evidently $ \lim_{N \to \infty} \rho(N) = 0.$ 
Let us also assume that only the non-trivial {\it infinite-dimensional case}
 will be considered, when an infinite multitude of Fourier coefficients 
$ f $ differs from zero, i.e. $ \forall N \ge 1 \ \Rightarrow \rho(N) > 0.$
 Otherwise our estimations will converge with speed $1/\sqrt{n}. $ \par
Now let us formulate the exact definition of an asymptotically optimal adaptive 
estimation (ADE), or, more precisely, a sequence of estimations. Let 
$ K(\theta), \ \theta \in \Theta $ 
be some set of Banach subspaces of space $L_2[0,1]; $ ( only the case when 
$ K(\theta) $ are compactly embedded in $ L_2[0,1] $ is non-trivial.) Set
$$
V(n,\theta) = \inf_{g(n)} \sup_{f \in K(\theta)} {\bf E} ||g(n) - f||^2,
$$
where $ \{g(n)\} $ is any sequence of estimations of $ f $. A sequence of 
adaptive estimations $ f(n) $ is called  {\it asymptotically optimal on 
the set of classes $ K(\theta) $ if} $ \forall \theta \in \Theta $
$$ 
\sup_n \sup_{f \in K(\theta)} {\bf E} ||f(n) - f||^2/V(n, \theta) < \infty.
$$
Of course, the quadratic function of losses $ l = l(g(n),f) = ||g(n) - f||^2 $ 
can be replaced by another loss function, non - negative, monotonically depending 
on the norm  $ ||g(n) - f||, $ so that
$$ 
\forall \epsilon > 0 \ \ \exists C(\epsilon) < \infty, \ \Rightarrow  \ l(z) \le 
C(\epsilon) \left(\exp(z^{\epsilon} \right) -1, \ z \ge 0.
$$
 Here is an important example of class $ K(\theta). $ Let $ \theta = 
\theta(N) $ be an arbitrary monotonically non-increasing numerical sequence 
such that $ \lim \theta(N) = 0; \ N \to \infty. $ Denote
$$ 
K(\theta) = \{f, \ ||f||^2(\theta) \stackrel{def}{=} \sup_{N \ge 1} 
\rho(f,N)/\theta(N) < \infty \}.
$$
 Relative to the norm $ ||\cdot||(\theta) $ class $ K(\theta) $ is a Banach 
space compactly embedded in $L_2[0,1], $ while the inverse is also true: any 
subspace compactly embedded in $L_2[0,1] $ is a subspace of some $K(\theta). $\par

 The value $\rho(f,N) $ is known and well studied in the {\it approximation 
theory. } Namely, $ \rho(f,N) = E^2_{N,2}(f),$ where $E_{N,p}(f) $ is the error 
of the best approximation of $ f $ by trigonometrical polynomials of power not 
exceeding $ N $ in $L_p $ metrics:  for $ g: [0,1] \to R^1 $ we shall denote 

$$ 
||g||_p = \left(\int_0^1 |g(x)|dx \right)^{1/p}, p \in [1, \infty);  \ 
||g||_{\infty} = \sup_{x \in [0,1]} |g(x)|, 
$$
 and closely connected with module of continuity of the form
$$ 
\omega_{p,2}(f^{(k)},\delta) = \sup_{h: |h|\le \delta}||f^{(k)}(x+h) - 2 f^{(k)}(x) 
+ f^{(k)}(x-h)||_p,
$$
(Timan A., 1960, p. 275);  arithmetical operations on the arguments of function 
$ f $ and their derivatives are understood modulo 1 (periodicity).\par
Everywhere below condition $ (\gamma1) $ will be considered fulfilled:

\begin{equation}
(\gamma1): \overline{\lim}_{N \to \infty} \rho(2N)/\rho(N) \stackrel{def}{=} 
\gamma < 1, \label{1}
\end{equation}
	
sometimes stronger conditions $ (\gamma) $ as well:
\begin{equation}
(\gamma): \exists \lim_{N \to \infty} \rho(2N)/\rho(N) \stackrel{def}{=} 
\gamma < 1;  \label{2}
\end{equation}

\begin{equation}
(\gamma0): \gamma = 0. \label{3}
\end{equation}

It is easy to show that from condition (1)  follows

\begin{equation}
\rho(N) \le C N^{-2\beta}, \ \ \  2 \beta \stackrel{def}{=} \log_2(1/\gamma) > 0.
\label{4}
\end{equation}

 In problem (R) it will be assumed that $\beta > 1/2. $ There are some grounds 
to suppose that at $\beta < 1/2 $ asymptotically optimal AE do not exist in 
the regression problem; for a similarly stated problem this was proved by 
Lepsky (Lepsky O., 1990). \par
 Also denote $ \kappa = \max(1,2\beta). $ Here and below the symbols $ C,C_r $ 
will denote positive finite constructive constants inessential in this context, 
$ \asymp $ is the usually  symbol, in detail:
$$
A(n) \asymp B(n) \ \Leftrightarrow C_1 \le \liminf_{n \to \infty} A(n)/B(n) \le 
$$
$$
\limsup_{n \to \infty} A(n)/B(n) \le C_2, \ \exists C_1, C_2 \in (0,\infty).
$$
 the symbol $ A \sim B $ means that in the given 
concrete passage to the limit $\lim A/B = 1.$ \par
{\bf Example 1.} Denote by $ W(C,\alpha,\beta) $ a class of functions $\{f\} $ 
such that
$$
\rho(f,N) \sim C N^{-2\beta} (\log N)^{\alpha}, \ \exists C,\beta >0; \alpha = 
const.
$$
Also denote $ W(\alpha,\beta) = \cup_{C>0}W(C,\alpha,\beta);$
$$ 
W(\beta) = W(0,\beta); \ \ W = \cup_{\beta > 0} W(\beta).
$$
 For the class of functions $ W $ condition $ (\gamma) $ is fulfilled. It is 
known from [19, pp. 275, 353] that $ f \in W(\alpha,\beta) $ if and only  if
 at $ \delta \to 0+, \ \delta \in (0, 1/2] $ 
$$
\omega_{2,2}(f^{[\beta]},\delta) \asymp \delta^{ \{\beta\}}|\log \delta|^{\alpha/2}, \ 
$$
$$
\forall j \le [\beta] \ \Rightarrow f^{(j)}(1-0) = f^{(j)}(0+0),
$$
where $ [\beta] $ denotes the integer part of $\beta $ and $ \{\beta\} = \beta -
[\beta]. $  At $ \{\beta\} = 0 $ function $ f^{[\beta]}(x) $ is assumed to be 
continuous.\par
{\bf Example 2.} Let us denote
$$
Z(\alpha,\beta) = \{f: \rho(f,N) \sim \alpha \ \beta^N \}, \alpha > 0, \ \beta 
\in (0,1);
$$
and $ Z = \cup_{\alpha > 0; \beta \in (0,1)} Z(\alpha,\beta) $. For functions of class $Z$ 
condition $(\gamma0) $ is fulfilled. Besides, functions of class $Z$ are 
analytical [20, p. 129].\par
 Denote for the problems $ {\bf R, D, S} $ respectively at $ j < n $ \ \ 
$ \hat{c}_j = (1/n) \times $

\begin{equation}
  (1/n) \sum_{i=1}^n y_j \varphi_j(x_i); \ \ \hat{c}_j = (1/n)
\sum_{i=1}^n \varphi_j(\xi_i); \ \ \hat{c}_{j} = \sum_{i=1}^{n-j}
\xi_i \xi_{i+j}/(n-j), \label{5}
\end{equation}
$ j = 0,2,4,\ldots $ and $ \hat{c}_j = 0 $ other case;
and for the regression problem
$$
c_j(n) = n^{-1} \sum_{i=1}^n \varphi_j(x_i), \ \ B(n,N)= \sum_{k=N+1}^{2N} 
c_k(n)^2 + \Delta_1 N/n;
$$
$ \Delta_1 = \sigma^2 = {\bf D} \xi_i; $  for problem D
$$
B(n,N) = \sum_{k=N+1}^{2N} c_k^2 + \Delta_2 N/n, \ \ \Delta_2 = 1;
$$
for the spectral problem
$$
B(n,N) = \sum_{k= N+1}^{2N} c^2_k + \Delta_3 N/n, \  \ \Delta_3 = ||f||^2;
$$
and again for all the problems set $ B(n) = $
$$
  \min_{N = 1,2,\ldots, [n/3]} B(n,N), \ N^0 = N^0(n) = 
\argm_{N = 1,2,\ldots,[n/3]} B(n,N);
$$
$$
A(n,N) = \rho(N) + \Delta_s N/n, \ A(n) = \min_{N = 1,2,\ldots,[n/3]} A(n,N),
$$
where $ s $ is the problem number.\par
 For instance, suppose that $ f \in W(C,\alpha,\beta) $, then $ A(n) \asymp 
n^{-2\beta/(2\beta+1)} (\log n)^{\alpha/(2\beta+1)} $, and in case 
$  f \in Z(\alpha,\beta) \ \Rightarrow A(n) \asymp \log n/n.$ \par
 Our notation should not be surprising, as it follows from the Bernstein 
theorem [20, p.242] and from condition $ (\gamma1) $
that all the introduced functionals $ \{B(n,N)\}, \ \{B(n)\} $ arising from 
different problems are mutually  $ \asymp $  equivalent. Besides, for the same 
reasons	
$$
A(n,N) \asymp B(n,N); \ A(n) \asymp B(n).
$$
 Let us make another additional assumption with regard to the class of 
estimated functions $ \{f\}: $
$$	
 (v): \forall v \ge 1, \  \forall N \in [1,N^0/v)]\cup[N^0 \cdot v,[n/3]] \ 
\Rightarrow \
$$

\begin{equation}
 B(n,N) - B(n) \ge C_1 (v-1)^2(1+C_2 |v-1|)^{-1} B(n). \label{6}
\end{equation}
(At $ v \ge N^0 $ the left interval of (6) is absent, at $ v \ge n/(3N^0) $ 
right interval of (6) is absent.) \par
 The classes of functions satisfying conditions $ (\gamma1) $ and $(v) $ will 
be called regular. Classes $ W $ and $ Z $ are regular.\par
 Apart from that it is clear that in the regression problem conditions must be 
imposed not only on the estimated function, but on the measurement errors 
$ \xi_i $ too. Two kinds of such conditions will be considered:
$$
(Rk): \ \exists k = 2,3,\ldots, \mu_{2k} \stackrel{def}{=} {\bf E} \xi_i^{2k} 
< \infty;
$$
(the power level) and the exponential level:
$$
(Rq): \exists q,Q \in (0,\infty), \Rightarrow {\bf P} (|\xi_i| > x) \le 
\exp \left(-(x/Q)^q \right), x > 0.
$$
The classical projective estimates by N. N. Tchentsov (Tchentsov N.N., 1972, p. 286)
 will be considered as estimates of the function $ f $:
$$
f(n,N,x) = \sum_{j=1}^N \hat{c}_j \varphi_j(x).
$$
Since, as shown by Tchentsov, 
$ {\bf E}||f(n,N,\cdot) - f(\cdot)||^2 \asymp B(n,N), $ the selection of the 
number of harmonics $ N $ optimal by order in the sense of 
$ L_2(\Omega) \times L_2[0,1] $ is given by the expression $ N = N^0(n) $
with the speed of convergence  $ f(n,N^0,\cdot) \to f(\cdot) $
in the above-mentioned sense $ \sqrt{A(n)}. $ I. A. Ibragimov and R. Z. 
Khasminsky (Ibragimov I., Khasminsky R., 1982)
 proved that no faster convergence exists on regular classes of 
functions given by the value $ \sqrt{A(n)}. $ \par
 However, the value $ \rho(f,N) $ or at least its order at $ N \to \infty $ are 
practically unknown as a rule. Below the adaptive estimation of $ f $ will be 
studied based only on observations $ \{\xi_i\} $ and using no apriory 
information regarding $ f $, and yet possessing the optimal speed of convergence 
 at apparently weak restrictions. Set 

\begin{equation}	
\tau(N) = \tau(n,N) \stackrel{def}{=} \sum_{k=N+1}^{2N} \hat{c}^2_k, \
 \ N(n) \stackrel{def}{=} \argm_{N \in (1, [n/3])} \tau(n,N),  \label{7}
\end{equation}

$$
\tau^*(n) = \min_{N \in (1,[n/3])} \tau(n,N), \ 
$$
{\it Our adaptive estimations $ \hat{f} $ in all considering problems  have a universal
 view:}
\begin{equation}
\hat{f} = f(n,N(n),x) = \sum_{j=1}^{N(n)} \hat{c}_j \varphi_j(x).  \label{8}
\end{equation}

 In case of a non-unique number of harmonics $ N(n) $ in (7) we choose the 
largest. Below the value $ N $ will always be in the set of integers numbers 
of segment $ 1,2,\ldots,[n/3]. $ \par
 Before proceeding to formulations and proofs let us clarify informally our 
idea for choosing $ N(n). $ It is easy to find by direct calculation that
\begin{equation}	
{\bf E} \tau(n,N) \asymp B(n,N), \ \ {\bf D} \tau(n,N) \asymp B(n,N)/n,  \label{9}
\end{equation}
and therefore 
$$ 
N \to \infty, N/n \to 0 \ \Rightarrow \sqrt{ {\bf D}\tau(n,N)}/{\bf E}\tau(n,N) 
\to 0.
$$
(In the case of the regression problem the condition $\beta > 1/2 $ is 
essential which is common in statistical research (Polyak B. at al., 1990, 1992),
(Lepsky O., 1990).  It 
follows from (9) that there are some grounds to assume
$$
\tau(n,N) \stackrel{a.s}{\asymp} {\bf E} \tau(n,N) \asymp A(n,N)
$$
and therefore
$$
N(n) = \argm_{N \le n/3} \tau(n,N) \sim \argm_{N \le n/3} {\bf E} \tau(n,N) 
= N^0(n).
$$
 Also note that the number of harmonics $ N(n) $ proposed by us is a random 
value (!) and that estimation (8) is non-linear by the totality of empirical 
Fourier coefficients $ \{\hat{c}_j\}. $ \par
 In the case of problem $ {\bf S } $ our estimation $ \hat{f} $ is homogeneos 
of degree 2 as a function of a initial data $ \{ \xi_i \} $ but also non - 
linear.\par 

\vspace{3mm}

\section{ Formulation of the main results.}

Let us denote
$$
{\bf P}_f(u) = {\bf P} \left(B^{-1}(n) ||\hat{f} - f||^2> u \right), \ u > C, \
 C \in (0,\infty).
$$

{\bf Theorem R.1(k).} {\it If the conditions $ (\gamma1) $ and $ \mu_4 < 
\infty $ is fulfilled in the regression problem, then}
$$
{\bf P}_f(u) \le C_1 \ \mu_4 \ u^{-1} \ \log^2(C_2 u), \ u > e/C_2,
$$

$$
C_1 = \min_{X \in (0,1/2)}(X(0.5 - X)^{-2} + X^{-1} ) \approx 7.221039 \ldots,
$$

$$
C_2 = \argm_{X \in (0,1/2)} (X(0.5 - X)^{-2} + X^{-1}) \approx 0.198340 \ldots.
$$
 This result was proved in (Bobrov P. at al., 1997),  but here the values of the 
constants have been corrected.\par
{\bf Theorem R.2(k).} {\it If the conditions $(\gamma1) $ and $(Rk) $ for some 
$ k = 3,4,\ldots $ is fulfilled in the regression problem, then }
$$
{\bf P}_f(u) \le 2^{2k} \ k^k \ \mu_{2k} \ u^{-k/2}, \ u > 0.
$$
{\bf Theorem R.3(q).} {\it In the conditions $(Rq), (\gamma1),(v) $ in the same 
problem at $ u > C = 2(1-\gamma)^{-1} Q $ the following inequality is true:}
$$
{\bf P}_f(u)\le 5 \exp \left[-C_1 \frac{ N^0(n) \ (u - C)/Q)^{q/(2q+4)} }
{ |\log B(n)| } \right].
$$

{\bf Theorem D.} {\it If in addition to the formulated conditions the 
boundedness of $ f $ is presumed, then in problem $ (D) $ at 
$ u \ge C = (1-\gamma)^{-1} $}
$$
{\bf P}_f(u) \le 5 \exp \left[ - C_2 \frac{ \sqrt{( u - C) N^0(n)}}
{|\log B(n)| } \right].
$$
{\bf Theorem S.} {\it If spectral density $ f(x) $ is bounded and conditions 
$ (\gamma1), (v) $ are fulfilled, then at $ u \ge C = (1-\gamma)^{-1} $ }
$$	
{\bf P}_f(u) \le 5 \exp \left[ - C_3 \frac{\sqrt{ (u - C) N^0(n)} }
{|\log B(n)|} \right].
$$
{\bf Theorem (R.k) \ a.s.} {\it If in problem R condition $(Rk) $ is fulfilled 
and the series }
$$
\sum_{n=1}^{\infty} n^{-k/2} A^{-k/2 + 2 \kappa}(n) < \infty,
$$ 
{\it converges, then}

\begin{equation}
\lim_{n \to \infty} \tau^*(n)/B(n) = 1, \label{10}
\end{equation}
{\it and if condition $(v) $ is also fulfilled, then}

\begin{equation}
\lim_{n \to \infty} N^0(n) /N(n)= 1. \label{11}
\end{equation}

(Recall that the convergence of a r.v. is understood in this paper only 
with probability 1).\newline

{\bf Theorem $ (Rq) \ a.s.$ } {\it If in the same problem under condition $ (Rq) $ 
for any $ \varepsilon > 0 $ the series}

\begin{equation}
\sum_{ \{n: \ 3 A(n)<1 \} } \exp \left(- \varepsilon \frac{(n A(n))^{q/(2q+4)}}
{|\log A(n)|} \right) < \infty, \label{12}
\end{equation}

{\it converges, then propositions (10) and (11) hold as well.} \par
{\bf Theorem (D)(S) a.s.} {\it Let for problems $ (D), (S), $ besides the 
above-formulated assumptions, condition (12) also be fulfilled with 
$ q/(2q+4) $ replaced by 1. Then the factual convergences of (10) and (11) 
are asserted here as well.}\\
(In comparison with (Bobrov P. at al., 1997)
 the exponent indices are significantly decreased.)\par

\vspace{3mm}
\section{ Auxiliary results.} 

\vspace{2mm}

 The technical apparatus for the proofs is the theory of so-called $ G(\psi) - $
spaces, i.e. Banach spaces of random values with rapidly diminishing tails of 
the distributions [16, 23]. For the reader's convenience the necessary 
information from that theory will be provided here without proof.\par
 A random value $ \eta $ determined, like all the other values in the present 
paper, on a fixed probability space, belongs to the space $ G(\psi) $, where 
$ \psi = \psi(m) $ is a function monotonically increasing on the set 
$ m \in (1, \infty) $ and finite at at least one value $ m > 1, $ if the norm 
$$
||\eta||(G(\psi)) \stackrel{def}{=} \sup_{k \ge 1}|\eta|_k/\psi(k) < \infty; \ 
\  \ |\eta|_k \stackrel{def}{=} {\bf E}^{1/k}|\eta|^k
$$
is finite. If $ \psi(m) = m^{1/q}, \ q = const > 0; $ then the corresponding 
space will be denoted $ G_p; \ p = q/(q-1); q = 1 \ \Rightarrow p = + \infty $
and the norm in it $ |||\eta|||_p; $ while $ \eta \in G_p $ then and only then, 
if
\begin{equation}
\exists C \in (0,\infty], \ \forall x > 0 \ \ {\bf P}(|\eta|>x) \le \exp 
\left(- C x^q \right). \label{13}
\end{equation}
 Now let $ \eta(t), \ t \in T, - $ be a separable random field, T an arbitrary 
set, and $ \sup_{t \in T} ||\eta(t)||_p \le 1.$ Introduce a so-called natural 
metric (more exactly semi-metric) $ d_p(t,s) = ||\eta(t)  - \eta(s)||_p $ and 
denote by $ {\bf N}(d_p,\varepsilon) $ the least number $d_p $ of spheres with 
radius $\varepsilon > 0 $ covering the entire set T. If the so-called entropic 
integral
$$
J = \int_0^1 \left(\log {\bf N}( d_p,\varepsilon)\right)^{1/q} d\varepsilon < \infty
$$
converges, then
\begin{equation}
||| \sup_{t \in T} |\eta(t)| \ |||_p \le C_1 + C_2 J.\label{14}
\end{equation}
 A similar result for spaces $L_k(\Omega) $ was obtained by G. Pizier (Pizier G.,
 1979 - 1980.) It is asserted that if
$$
1) \exists k > 1 \ \Rightarrow \sup_{t \in T} {\bf E}|\eta(t)|^k \le 1;
$$
$$
2) I \stackrel{def}{=} \int_0^1 {\bf N}^{1/k}(r_k,\varepsilon) 
d \varepsilon < \infty,
$$
where $ r_k(t,s) = |\eta(t) - \eta(s)|_k,$ then
\begin{equation}
| \ \sup_{t \in T}  |\eta(t)| \ |_k \le C_1 + C_2 I. \label{15}
\end{equation}

\vspace{2mm}
	
\section{ Proofs }
\vspace{2mm}

 The proofs of the theorems referring to different problems are similar. The 
assertions referring to problem $ R,$ which is the most complicated, will be 
proved below in detail, and after that the changes will be indicated that arise 
in considering problems $ D $ and $ S.$ Some additional notations have to be 
introduced: for $ f: [0,1] \to R^1 $  and $ p \ge 2 $ we shall denote
$$
 ||f||_{p,d} = \left\{n^{-1} \sum_{i=1}^n |f(x_i)|^p \right\}^{1/p},
$$
while in the case of $p=2 $ the index $ p $ of the norm sign will be omitted. 
Further,
$$
\Phi(N,f,x) = \Phi(N,x) = \Phi(N) = \sum_{j=1}^N c_j \varphi_i(x) 
$$
are partial Fourier sums for the function $ f(x) $,
$$
T(N) = T(N,x) = \Phi(2N,f,x) - \Phi(N,x), \ N \le n/3.
$$
{\bf Lemma 1.} For all $ p \ge 2 $
$$
||T(N)||_p \asymp ||T(N)||_{p,d} \le C N^{1/2 - 1/p} \ \ \sqrt{\rho(f,N)}.
$$
{\bf Proof.} The first assertion follows from the fact that $ N \le n/3 $ and 
from the Bernstein inequality [19, p. 245]. The other uses the Nikolsky 
inequality (Timan A., 1960, p. 245):
$$
n^{-1} \sum_{i=1}^n |T(N,x_i)|^p = ||\Phi(2N) - \Phi(N)||^p_{p,d} \le 
2^p ||\Phi(2N) - \Phi(N)||^p_p \le 
$$

$$
 6^p N^{p/2 - 1} ||\Phi(2N) - \Phi(N)||^p_p \le 6^p  N^{p/2 - 1}
(\rho(N) - \rho(2N))^{p/2} <
$$

$$
 6^p N^{p/2 - 1} \rho^{p/2}(N).
$$

{\bf Lemma 2.} Let us consider on the set $ S = \{1,2,\ldots, n\} $ the metric
$$
d(N_1,N_2)= |\rho(N_1) - \rho(N_2)| + n^{-1}|N_1 - N_2|.
$$
 It is asserted that the entropy of the set $ S $ in the metric $ d,$ i.e.
$ H(S,d,\varepsilon) = \log {\bf N}(S,d,\varepsilon), \ \varepsilon \in (0,1] $
satisfies the inequality
$$
H(S,d,\varepsilon) \le C + \kappa |\log \varepsilon|, \ \kappa = \max(1,2\beta).
$$
{\bf Proof.} Set $ K = C \left[\varepsilon^{-1/\beta}\right] $ and consider 
$ S(\varepsilon) $ - the 
net $ S $ in the metric $ d $ of the form
$$
S(\varepsilon) = \left[ \left( \{1,2,\ldots,K\} \right) \cup  
\left(\cup_j\{[nj\varepsilon/2]\}\right) \right] \cap S.
$$
Calculation of the number of elements in $ S(\varepsilon) $ convinces us of 
the correctness of the lemma.\par
 The central moment in all the further considerations is the so-called 
{\it expansion of the basic functional} $\tau(n,N). $ In all the three problems 
under consideration $\tau(n,N) $ is of the form $ \tau(n,N) = $
$$
= {\bf E}\tau(n,N) +2 \Psi_1(N) + \Psi_2(N), \ \ {\bf E} \tau(n,N)\sim B(n,N),
$$
$$
{\bf E} \Psi(N) = {\bf E} \Psi_2(N) = 0; \ \Psi_s(N) = \Psi_s(n,N),
$$
where in the case of the regression problem $ \Psi_1(N) = n^{-1} \times $
$$
\times \sum_{i=1}^n \xi_i T(N,x_i), \ \Psi_2(N) = n^{-2} \sum \sum_{i,j=1}^n
a_{i,j}(n,N)(\xi_i\xi_j - {\bf E} \xi_i \xi_j),
$$
$$
a_{i,j}(n,N) = D_{2N}(x_i,x_j) - D_N(x_i,x_j),
$$
$ D_N(x,y) = \sum_{j=1}^N \varphi_j(x) \varphi_j(y) - $ is the Dirichlet kernel.\par
 It is easy to obtain by direct calculation for problem R (and then for the 
remaining problems) that
$$
{\bf D} \Psi_1(N) \asymp \rho(N)/n, \ \ \ {\bf D} \Psi_2(N) \asymp N/n.
$$
{\bf Lemma 3.} In the regression problem under conditions $ (Rk) $ the 
following inequality holds:
\begin{equation}
|\Psi_1(N)|_{2k} \le C \  k \ \mu_{2k}^{1/2k} \ \sqrt{\rho(N)/n}, \ k = 2,3,\ldots
. \label{16}
\end{equation}
{\bf Proof.} We shall apply the moments inequalities for the sums of centered 
independent variables  $ \{\varepsilon_i \}, \ i = 1,2,\ldots,n $ at (Rosental H., 1970),
(Johnson W.B., Shechtman G, Zinn J. at al., 1985): $ p \ge 2 \ \Rightarrow $
$$
|\sum \varepsilon_i|_{p} \le 3 (p/\log p) \max \left( |\sum \varepsilon_i|_2, 
(\sum |\varepsilon_i|_p^p)^{1/p} \right).
$$
 Here $ \varepsilon_i = \xi_i \ T(N,x_i), \ \sum_i = \sum_{i=1}^n, \
\ p = 2k. $ As long as, on the basis lemma 1,
$$
\sum_i |T|^p(N,x_i) = n ||T(N)||_{p,d}^p \le C n ||T(N)||_p^p \le 
$$
$$ 
 C n N^{p/2 - 1} ||T(N)||_2^{p} = C n N^{p/2 - 1} \rho^{p/2}(N),
$$
we obtain the conclusion of lemma 3.\par

{\bf Lemma 4.} In the same problem and in the same assumptions
\begin{equation}
|\Psi_2(N)|_k \le C \ k \ \mu^{1/k}_{2k} \ \sqrt{N}/n.\label{17}
\end{equation}
{\bf Proof.} It is sufficient to prove (17) for even $ k $, while for the odd 
ones it is necessary to consider the moment of order $ k+1 $ and make use of 
the Lyapunov inequality. The functional $ \Psi_2(N) $ is the quadratic centered 
form from the random values $ \{\xi_i\}, \ i = 1,2,\ldots,n. $ In order to 
estimate its k-th moment we will estimate its cumulant of the same order. 
According to [29, p. 101],
$$
\frac{  \Gamma_k(\Psi_2(N))}{ {\bf D}^{k/2} \Psi_2(N)} \le C^k k^k \mu_{2k} \cdot
\left( \frac{W_n}{ {\bf D}^{(k-2)/2} \Psi_2(N)} \right)^{k-2}, \ k \ge 4,
$$
where $ \Gamma_k(\xi) $ denotes the k-th semi-invariant of the value $ \xi, $
$ W_n = $
$$
 \max_i \sum_{j=1}^n |a_{i,j}(n,N)| \le \max_{y \in [0,1]} \sum_{j=1}^n 
\left|\sum_{l=N+1}^{2N} \varphi_l(y) \varphi_l(x_i) \right|.
$$
Analogously to the estimations of the Lebesgue constants in the theory of 
trigonometrical series we can estimate $ W_n \le C \log N/n, $ and 
consequently at $ k \ge 4 $
$$
{\bf \Gamma}_k(\Psi_2(N)) \le C^k \ k^k \ \mu_{2k} \ N^{-1} \ \log N \cdot {\bf D}^{k/2}
(\Psi_2(N)).
$$
Proceeding by the well-known Leonov - Shiryaev formulas (Shiryaev A.N., 1989, p.311)
 from semi-invariants to moments, we obtain the proposition of lemma 4.\newline
{\bf Lemma 5.} Under the conditions of Lemmas 3 and 4 the following 
inequalities hold respectively:
$$
|\tau(n,N) - {\bf E} \tau(n,N)|_k \le C \  k \ \mu_{2k}^{1/k} \ \sqrt{ A(n,N)/n},
$$
and if condition $ (Rq) $ is fulfilled, then on the basis of the properties of 
spaces $  G(\psi) $
$$
||| \tau(n,N) - {\bf E} \tau(n,N)|||_r \le C \sqrt{A(n,N)/n}, \ r = q/(q+2).
$$
 The assertion of the lemma  5 it follows from the inequality of the triangle for the used
norms. \par
 Let us  consider the centered and normalized random  field
$$
\zeta(N)= \zeta(n,N) = \sqrt{n/(A(n,N)} \ [\tau(n,N) - {\bf E}\tau(n,N)],
$$
so that $ {\bf E}\zeta(N)=0, \ \sup_{N \le [n/3]}  || \zeta(N) ||_r \le C. $ \par
{\bf Lemma 6.}
$$
(Rk) \Rightarrow |\max_{N} |\zeta(N) \ |_k \le C A^{-2\kappa/k}(n), \ 
k> 3 \kappa;
$$
$$
(Rq) \Rightarrow ||| \max_N |\zeta(N)| \ \ |||_r \le C |\log A(n)|^{1/r}.
$$
 The {\bf proof} will be given for the second assertion alone, as the first is 
simpler because the spaces $ L_k(\Omega) $ are more customary. We obtain on the basis of 
lemma 5, put
$$ 
\Psi(N) = 2\Psi_1(N) +\Psi_2(N) = \tau(n,N) - {\bf E}\tau(n,N):
$$

$$
n^{-1/2} \left(\zeta(N_1) - \zeta(N_2) \right) = 
\frac{ \left( \sqrt{A(n,N_1)}-\sqrt{A(n,N_1)}\right) \Psi(N_1)}
{\sqrt{A(n,N_1) A(n,N_2)} }+
$$

$$
 [\Psi(N_1)-\Psi(N_2)] / \sqrt{A(n,N_2)} \stackrel{def}{=} \zeta_1 + \zeta_2;
$$

$$
||| \zeta_2 |||_r \le C \sqrt{ |A(n,N_1)- A(n,N_2)|}/(A(n)\sqrt{n}) \le 
$$

$$
 C\sqrt{ |\rho(N_1) - \rho(N_2)| + n^{-1}|N_1 - N_2|}/(A(n)\sqrt{n}),
$$
since $ A(n,N) \ge A(n), a \ge b \ge 0 \Rightarrow \sqrt{a} - \sqrt{b} \le 
\sqrt{a-b},$
$$
||| \Psi(N_1) - \Psi(N_2) |||_r \le C \sqrt{|A(n,N_1)- A(n,N_2)|/n}.
$$ 

Further, $ |||\zeta_1 |||_r \le ||| \Psi(N_1) |||_r \times $
$$
 \left[|\sqrt{A(n,N_1)} - \sqrt{A(n,N_2)}| \right] \cdot
\left[\sqrt{A(n,N_1) A(n,N_2)} \right]^{-1/2} \le 
$$
$$
 C \sqrt{|\rho(N_1) - \rho(N_2)| + n^{-1}|N_1 - N_2|}/(A(n)\sqrt{n}),
$$
 since $ ||| \Psi(N_1) |||_r \le C /\sqrt{n}. $ The random field $ \zeta(N) $ 
is thus bounded in the norm $ ||| \cdot |||_r $ and
$$
d_1(N_1,N_2) \stackrel{def}{=} ||| \zeta(N_1)-\zeta(N_2) |||_r \le C 
\sqrt{d(N_1), N_2)}/A(n).
$$
Since $ H(S,d_1,\varepsilon) \le H(S,\sqrt{d}/(CA(n)), \varepsilon) = $
$$
 H(S,d,(C \varepsilon A(n))^2) \le C_1 + 2\kappa |\log \varepsilon| + 
 2 \kappa |\log A(n)|.
$$
the assertion of the lemma follows from the properties of the spaces $ G(\psi) $ 
(13, 14, 15).\par
 Inequalities (16) and (17) can be reformulated as follows in forms more 
convenient for further application. Under conditions $ (16) $ and $ (17) $ 
respectively the sequence $ \{\tau(n,N)\} $ can be expended into
$$
\tau(n,N) = {\bf E}\tau(n,N)+ \sqrt{{\bf E}\tau(n,N)/n} \ k \ \mu_{2k}^{1/k} \
(A(n))^{-2\kappa/k} \ \nu(n,N),
$$
 where ${\bf E} \nu(n,N) = 0,$ 
\begin{equation}
\sup_n {\bf E} \max_N |\nu(n,N)|^k = C < \infty; \label{18}
\end{equation}
and in the other case $ (Rq) $
$$
\tau(n,N)={\bf E} \tau(n,N)+ \sqrt{{\bf E} \tau(n,N)/n}\cdot|\log A(n)|^{1/r}
\cdot \nu(n,N),
$$
$$
{\bf E} \nu(n,N) = 0;  \ \ \sup_n ||| \max_N |\nu(n,N)| \ \ |||_r = C < \infty.
$$
{\bf Lemma 7.} Let $ M $ be some subset of an integral segment 
$ S = [1,2,\ldots,n], $ $ \overline{M} = S \setminus M, \ \pi(M)= {\bf P}
(N(n) \in M), $
$$
v =v(n,M) \stackrel{def}{=} \inf_{N \in M} B(n,N)/B(n) \ge 2.
$$
Then under conditions $ (Rq) $
\begin{equation}
\pi(M) \le 2\exp \left(- C \  (v n B(n))^{r/2} /
 |\log B(n)|  \right) \label{19}
\end{equation}
and under conditions $ (Rk) $
\begin{equation}
\pi(M) \le \frac{C^k \ k^k \ \mu_{2k}^k}{ v^{k/2} \ n^{k/2} \ B^{k/2 - 2 \kappa}(n)}.
\label{20}
\end{equation}
{\bf Proof.} We obtain for the case of $ (Rq) $, denoting $ \overline{\nu} = 
\max_{N \in S} |\nu(n,N)|: $  
$$
\pi(M) = {\bf P}(N(n) \in M) = {\bf P} (\min_{N \in \overline{M}} \tau(n,N) >
\min_{N \in M} \tau(n,N)) =
$$

$$
 {\bf P} \left( \min_{N \in \overline{M}} ( B(n,N) + \sqrt{B(n,N)/n}
 \ |\log B(n)|^{1/r} \ \nu(n,N)) \right) >
$$

$$
> \min_{N \in M} \left( B(n,N) + \sqrt{B(n,N)/n} \ |\log B(n)|^{1/r} \ \nu(n,N) 
\right)  \le 
$$

$$
 {\bf P} ( B(n) + \sqrt{B(n)/n} \ |\log B(n)|^{1/r} \ \overline{\nu} >
$$

$$
> v B(n) - \sqrt{v B(n)/n} \ (|\log B(n)|)^{1/r} \ \overline{\nu} ).
$$
 We find  solving the inequality under the probability symbol relative to 
$ \overline{\nu} $ (the case of $ v \ge C /B(n) $ is trivial):
$ \pi(M) \le $
$$
  {\bf P}(C \ \overline{\nu} \ (1 + \sqrt{v}) \sqrt{B(n)/n} \ |\log B(n)|^{1/r} 
\ge (v-1) B(n) ) \le 
$$

$$
 {\bf P}\left( \overline{\nu} \ge C \frac{v-1}{ \sqrt{v}+ 1}  \ 
\frac{ \sqrt{n B(n)} } {|\log B(n) |^{1/r}} \right) \le 
 {\bf P} \left( \overline{\nu} \ge C \sqrt{v} \frac{\sqrt{n B(n)}}
{|\log B(n)|^{1/r} } \right).
$$
 Using the estimations of lemma 6, we arrive at inequalities (19) and (20). 
The case of $ (RK) $ is considered analogously.\par
{\bf Proof of Theorem (R.k) \ a.s.} It follows from expansion (18) that
$$
\tau^*(n) \le B(n) + \sqrt{B(n)/n} \ C(k) \ B^{-2\kappa/k}(n) \ \overline{\nu},
$$
therefore 
$$
 \left|\frac{\tau^*(n)}{B(n)} - 1 \right| \le \frac{C(k)}{ \sqrt{v}
B^{1/2 - 2 \kappa/k}(n)  \overline{\nu}}.
$$
We receive in according to the Chebyshev inequality:
$$
{\bf P}_n(\varepsilon) \stackrel{def}{=} {\bf P} \left( \left| \frac{\tau^*(n)}
{B(n)} - 1 \right| > \varepsilon \right) \le \frac{C(k)}{\varepsilon^k n^{k/2}
B^{k/2 - 2 \kappa}(n)}.
$$
Since for any $ \varepsilon > 0 $ the series $ \sum_n {\bf P}_n(\varepsilon) $ 
converges, the first assertion to be proved follows from the Borel - Cantelli 
lemma. The other is proved analogously if it is taken into account that 
$ N \ge N^0(n) (1 + \varepsilon), \ \varepsilon \in (0,1] $ and condition 
$ (v) $ lead to the inequality $ B(n,N) \ge (1 + C \varepsilon^2) B(n), $ 
$ \varepsilon \in (0,1) $ and lemma 7 is applied.\par
 Analogously we can prove the theorem {\bf (Rq)\ `a.s}, on the basis of 
inequality: 
$$
{\bf P}_n (\varepsilon) \le \exp \left( - C \varepsilon^{r} \ (n A(n))
/ |\log A(n)| \right).
$$
 {\bf Remark 1.} Let us note, and use it below, a slight difference in the 
behaviors of the values $ \tau(n,N) $ and $ N(n) $ which consists in the 
peculiarity of condition $ (v). $ At $ v \ge 2 $ we have (under the same 
conditions $(Rq),(v): $
$$
\max \left( {\bf P} \left(\frac{N(n)}{N^0(n)} \le \frac{1}{v} \right),
{\bf P} \left( \frac{N(n)}{N^0(n)} > v \right) \right) \le \exp \left( - C v^r
\frac{(n A(n))^{r/2} }{|\log A(n)| } \right).
$$
 An analogous estimation for the probability $ {\bf P}(\tau^*(n)/B(n) > v) $ 
holds even without condition $ (v).$ \par
{\bf Remark 2.} The {\it consistency } of the proposed estimations in the 
above-mentioned sense follows from the assertions already proved. Indeed, 
since 
$$ 
A(n) \le A(n,[\sqrt{n}]) \le C n^{-1/2} + \rho([\sqrt{n}]) \to 0, 
$$ 
 then $ N^0(n) \to \infty, \ N^0(n) / n \to 0, $ because otherwise the value 
$$ 
A(n) = A(n,N^0(n)) \asymp N^0(n) / n+ \rho(N^0)
$$ 
would not tend to zero. \par
Since $ N(n)/N^0(n) \to 1, $ then $ N(n) \to \infty $ and analogously 
$ N(n)/n \to 0, $ which proves the consistency of $ \hat{f}.$ \par
{\bf Proof of Theorem R.3(q).} (The previous theorem is proved analogously). 
Note that because of condition $ (\gamma1) $ for $ n \ge n_0 > 2 $
$$
||\hat{f} - f||^2/B(n) = B(n,N(n))/B(n) + \Psi_3(N(n)) /B(n) \le 
$$
$$
(1-\gamma)^{-1} \tau^*(n)/B(n) + \Psi_3(N(n))/B(n) =
$$
$$
 (1-\gamma)^{-1}  + (\tau^*(n) / B(n) - 1) + \Psi_3(N(n)) /B(n),
$$
where, as can easily be seen, $ \Psi_3(N) =  $
$$
\sum_{l,s = 1, l \ne s}^n V(\xi_l,\xi_s), \ V(x,y) = \sum_{j=1}^N
(\varphi_j(x) - c_j)(\varphi_j(y) - c_j),
$$
and has the same form and the same estimation as $\Psi_2(N). $ \par
Then we will use the elementary inequality 
$ {\bf P}({\bf A}) \le {\bf P}( {\bf ABC}) + {\bf P}({\bf \overline{B}}) +
{\bf P}({\bf \overline{C}}), $ in which $ {\bf A,B,C} $ are events. Setting
$ {\bf A} = $
$$
 \{ ||\hat{f} - f||^2/B(n) > u\}, \ {\bf B} = \{ 1/v \le \tau^*(n) / 
B(n) \le v \},
$$
$ {\bf C} = \{ N^0 /v \le N(n)\le v N^0(n),  \} $ we have  at $ v \in (2,u - C): $ 
$$
{\bf P}_0 \stackrel{def}{=}{\bf P} ({\bf ABC} ) \le {\bf P}(v+ 
\max_{N \le v N^0} |\Psi_3(N)|/B(n) > u).
$$
We find analogously to lemma 4: $ |||\Psi_3(N)|||_r \le $
$$
 C \sqrt{N}/n, \ \ |||\Psi_3(N_1) - \Psi_3(N_2) |||_r \le \sqrt{|N_1 - N_2|}/n,
$$
and since the entropy integral converges, then (see (15))
$$
||| \max_{N \le v N^0(n)} |\Psi_3(N)| \ \ |||_r \le C \sqrt{v N^0(n)}/n \le 
C \sqrt{v} /\sqrt{N^0(n)}.
$$
 Using triangular inequality for the $ G(\psi) \ - $ norms  we obtain:
$$
{\bf P} \left(\left| \frac{\tau^*(n)}{B(n)}-1 \right| + 
\left| \frac{\Psi_3(N)}{B(n)} \right| > v \right) \le
 \exp \left( - C_5 \ v^r \ \frac{(n A(n))^{r/2}}{ |\log A(n)|} \right).
$$
We obtain therefore, based on the properties of spaces $ G(\psi): $
$$
{\bf P}_0 \le \exp \left( - C_6 \ \left( v^{-1/2} \ ( u - C -v) \ \sqrt{N^0(n)}  
\right)^r \right).
$$
 The other probabilities $ {\bf P}({\bf {\overline{B}}}),
{\bf P}({\bf {\overline{C}}}) $ were estimated above, and we find by summing
$ ( C = 1/(1-\gamma) \ ): $
$$
{\bf P}_f(u) = {\bf P} ({\bf {A}}) \le \exp \left( - \left(C_1(u- C - v)
\sqrt{ N^0(n) / v} \right)^r  \right) + 
$$

$$
+ \ 4 \ \exp \left( - C_2 v^{r/2} (n A(n))^{r/2} \right). \  
$$
 Taking into account that $ n B(n) > N^0 $ and choosing \\
$ v = C_4 (u - C), \ C_4 \in (0,1), $ we arrive at the assertion of the 
theorem.\par
{\bf We proceed now to the problem of estimating density (D).} The functional 
$ \Psi(n,N) = \Psi(n,N) $ has in it the following form:
$$
\Psi(N) = n^{-1} \sum_{i=1}^n \sum_{j=N+1}^{2N} (c_j \varphi_j(\xi_i)-
c_j^2).
$$
Using the Rosental inequality once more, we obtain
$$
{\bf E}(\Psi_1(N))^{2k} \le 2 \ C(2k) \ n^{-k} \ {\bf E} \left(\sum_{j=N+1}^{2N}
c_j \varphi_j(\xi_1) \right)^{2k} =
$$
$$
2 \ n^{-k} \ C(2k) \int_0^1  \left(\sum_{j=N+1}^{2N} c_j \varphi_j(x) \right)^{2k}
f(x) \ dx \le 
$$
$$
 C\cdot C(2k)n^{-k} \int_0^1 \left(\sum_{j=1}^{2N} c_j \varphi_j(x) 
\right)^{2k} dx,
$$
since $ f(x) $ is presumed to be bounded. Then, since
$$
|| \sum_{j=N+1}^{2N} c_j \varphi_j(x) || = || \Phi(2N,x) - \Phi(N,x)|| \to 0, 
\ N\to \infty,
$$
we have in according to the Riesz theorem (Timan A., 1960, p. 305)
$$
||\Phi(2N) - \Phi(N)||^{2k}_{2k} \le C^k \ k^k \ ||\Phi(2N) - \Phi(N)||^{2k} =
$$

$$
  C^k \ k^k \ (\rho(N) - \rho(2N))^k < C^k \ k^k \ \rho^k(N),
$$
so that
\begin{equation}
{\bf E}\Psi^{2k}_1(N) \le C^k \ k^k \ n^{-k} \ \rho^k(N).\label{21}
\end{equation}

 In the language of $ G(\psi) \ $ - spaces inequality (21) means that
$$
||| \Psi_1(N) |||_r \le C \rho(N)/n.
$$
It is proved analogously that
$$
||| \Psi_1(N_1) - \Psi_1(N_2) |||_r \le C |\rho(N_1) - \rho(N_2)|/n.
$$
The functional $ \Psi_2(N) = \Psi_2(n,N) $ has the form
$$
\Psi_2(N) = \sum_{l=1}^n \sum_{s=1}^n U(\xi_l,\xi_s),
$$
where
$$
U(x,y) = U(N,x,y) = \sum_{j = N+1}^{2N} (\varphi_j(x) - c_j)(\varphi_j(y)-c_j),
$$
and is consequently a so-called $ U- $ statistic with the kernel 
$ U = U(N,x,y). $ At the same time our $ U- $ statistic is singular. 
The asymptotics of the moments of this kind of statistics and the limiting 
distribution for them are to be found e.g. in (Korolyuk V.S. at al., 1989),
(Ronzin A., 1982).  However, here we 
need non-asymptotic estimations from above, and therefore additional reasoning 
will be required. Note first of all that
\begin{equation}
{\bf E}|U(\xi_1,\xi_2)|^m \le C^m N^{m-1}, \ m = 3,4, \ldots.\label{22}
\end{equation}
Let us prove (22).
$$
{\bf E} |U(\xi_1,\xi_2)|^m \le4^m C^m \int_0^1 \int_0^1 |D_{2N}(x,y)|^m 
f(x) f(y) dx dy \le 
$$
$$
\le C^m \int_0^1 \int_0^1 |D_{2N}(x,y)|^m dx dy,
$$
where, let us recall, $ f $ is bounded and $ D_N $ is the Dirichlet kernel.
The last integral is easily estimated and we arrive at (22). Then on account 
of the singularity of the statistics we have in the case of even $ k $:
$$ n^k{\bf E}\Psi_2^k(N) \le 
 \sum_{i_1=1}^n \sum_{j_1=1}^n \ldots \sum_{i_k=1}^n
\sum_{j_k=1}^n {\bf E} U(\xi_{i_1},\xi_{j_1}) \ldots U(\xi_{i_k}\xi_{j_k}) \le 
$$

$$
 C^k k^k \sum_{i_1 = 1}^n \sum_{j_1=1}^n\ldots \sum_{i_{k/2}=1}^n \sum_{j_{k/2}=1}^n
{\bf E} U^2(\xi_{i_1},\xi_{j_1}) \ldots {\bf E} u^2(\xi_{i_{k/2}},\xi_{j_{k/2}}) \le
$$

$$
 C^k k^k |U(\xi_1,\xi_2)|^k_k \le C^k k^k N^{k/2}.
$$
In the case of odd $ k $ we consider the moment of order $ k+1; $ in 
[16, p. 42] the equivalence is proved of the norms $ G(\psi), $ constructed by 
even moments alone, to the initial norm.\par
 Thus $ ||| \Psi_2(N) |||_r \le C \sqrt{N}/n, $ and the further course of 
reasoning is fully analogous to the ground for estimation of regression.\par 

{\bf Consider now the problem of spectral statistics (S).} It turns out 
unexpectedly that the reasoning here is even simpler than in problem (D). The 
fact is that the initial sequence $ \{\xi_i\} $ is assumed to be Gaussian, the 
empirical Fourier coefficients $  \hat{c}_k $, i.e. empirical correlation 
coefficients, are quadratic functionals from the trajectory $ \{\xi_i \}, \
\ i = 1,2,\ldots,n $, while the functional $ \tau(n,N) $ is a polynomial 
functional of the 4th power and therefore has the expansion
$$
\tau(n,N) = {\bf E} \tau(n,N) + \sum_{m=1}^4 \Psi_m(n,N),
$$
where the expansion components are not correlated between themselves and 
$ \Psi_m $ can be written as an $ m \ - $ dimensional stochastic Ito-Wiener 
integral according to the orthogonal Gaussian measure. At the same time
$$
C N /n \ge {\bf D} \tau(n,N) = \sum_{m=1}^4 \Psi_m(n,N),
$$
therefore $ {\bf D} \Psi_m(n,N) \le C N/n. $ The Plikusas theorem 
(Plikusas A., 1981)  asserts 
that the distribution $ \Psi_m(n,N)  $ is estimated only through dispersion: 
$$
||| \Psi_m(n,N) |||_{2/m} \le C(m) {\bf D}^{1/2} \Psi_m(n,N) \le C \sqrt{N/n},
$$
consequently $ ||| \tau(n,N) - {\bf E} \tau(n,N) |||_{1/2} \le C \sqrt{N/n} $. 
Analogously considering the dispersion of the value
$$
\zeta(N) = \zeta(n,N) = \sqrt{n/A(n,N)} \  [\tau(n,N) - {\bf E} \tau(n,N)],
$$
we find that $ {\bf D} \zeta(N) \le C /B(n) $ and therefore 
$ |||\zeta(N)|||_{1/2} \le C/B(n), $ and the difference 
$ \zeta(N_1) - \zeta(N_2) $ is estimated likewise:
$$
||| \zeta(N_1) - \zeta(N_2) |||_{1/2} \le C \sqrt{d(N_1,N_2)} /B(n).
$$
 As a result we obtain for the functional $ \tau(n,N) $ expansion (18), which 
is of key importance for us:
$$
\tau(n,N) = {\bf E} \tau(n,N) + \sqrt{{\bf E} \tau(n,N)/n} \cdot \log^2 B(n) 
\cdot \overline{\nu},
$$

$$
\sup_n ||| \max_{N \le n/3} | \overline{\nu} | \ \ |||_{1/2} = C < \infty.
$$
The other details of the proof are analogous to the case of regression and 
ought to be omitted.\par

\vspace{2mm}

\section{ Adaptive confidence intervals.} 
\vspace{2mm}
Let us now describe the use of our 
results for the construction of ACI. Note first of all that the probability 
$ {\bf P}_f(u) $ with rather weak conditions (except $ (Rk) $ ) in all the 
considered problems permits estimation of the form
\begin{equation}
{\bf P}_f(u) \le 5 \exp \left(- \varphi(C,N^0(n),B(n)) u^{r/2}) \right), 
\ u > C. 
\label{23}
\end{equation}
As proved above, the values $ N^0, B(n) $ have respective consistent estimates
$$
N^0(n) \approx \argm_{N \le n/3} \tau(n,N), \ B(n) \approx \min_{N \le n/3}
\tau(n,N) = \tau^*(n).
$$
The value $ C $ also depends on $ \gamma $ and on the constants $ C_j $ 
appearing in the definition of condition $ (v). $ With very weak conditions 
they can also be estimated consistently by the sampling in the following way. 
Set $ M = M(n) = \left[\exp(\sqrt{\log n}) \right] $; then, if conditions 
$ (\gamma), (v) $ are fulfilled, a system of asymptotic equalities can be 
written:
$$
\tau(M) - \Delta_s M/n \sim (1-\gamma)\rho(M);
$$
$$
\tau(2M) - 2 \Delta_s M/n  \sim \gamma (1-\gamma) \rho(M);
$$
$$
\tau(4M) - 4 \Delta_s M/n \sim \gamma^2 (1-\gamma) \rho(M).
$$
Solving this system, we find the consistent $ (mod \ \ {\bf P}) $ estimate of 
$ \gamma: $
$$
\hat{\gamma} = \frac{\tau(4M) - 2 \tau(2M)}{\tau(2M) - 2 \tau(M)}.
$$
(The parameter $ \Delta_s $ can also be estimated consistently, but that is not 
necessary for us). Further, since
\begin{equation}
\frac{\tau(n,N(n)(1+v))}{\tau^*(n)} \sim \frac{B(n,N(n)(1+v))}{B(n)}
\sim \frac{C_1 v^2 }{1+C_2 v },\label{24}
\end{equation}
the constants $ C_1, C_2 $ can be determined from (24), for instance by the 
least-squares method. Substituting the obtained estimates of all the parameters 
into (23), we get the estimate of the confidence probability
\begin{equation}
\ \ {\bf P}_f(u) \le 5 \exp \left( - \phi(C(\hat{\gamma}, \hat{C}_1,\hat{C}_2), N(n),\tau^*(n)) u^{r/2} \right)
\stackrel{def}{=} \hat{ {\bf P}}_f(u).\label{25}
\end{equation}
then, equating the right-hand part (25) of the unreliability of the confidence 
interval $ \delta $ to, say, the magnitude 0.05 or 0.01, we calculate 
$ u = u(\delta) $ from the relation
$$
\hat{ {\bf P}}_f(u(\delta)) =\delta
$$
and obtain approximately the {\it adaptive confidence interval} for $ f $ 
reliability $ 1 - \delta $ of the form
\begin{equation}
|| \hat{f} - f||^2 \le u(\delta) \min_{N \le n/3} \tau(n,N).\label{26}
\end{equation}
But for a rough estimate of the error from replacing $ f $ by $ \hat{f} $ the 
following quite simple method can be recommended. Since
\begin{equation}
\frac{||\hat{f} - f||^2}{ B(n)} = \frac{A(n,N(n))}{B(n)} + \frac{\Psi_3(N(n))}
{B(n)},\label{27}
\end{equation}
and the second term in the right-hand part of (27) a.s. tends to 
zero, while the first term, if conditions $ (\gamma), (v) $ are fulfilled, has 
$ 1/(1 - \gamma) $ as its limit, we thus prove the following assertion 
apparently well known to specialists in nonparametric statistics for 
non-adaptive estimation: \newline
{\bf Theorem c.i.} {\it If the following conditions are fulfilled in our 
problems: $ $ in problem $ R \ (Rq), (\gamma), (v) $ or $ (\gamma), (v) $ in 
problems $ D, \ S, $ then}
\begin{equation}
\lim_{n \to \infty} ||\hat{f} - f||^2 /B(n) < 1/(1 - \gamma).\label{28}
\end{equation}
In order to construct an adaptive confidence interval assertion (28) can be 
reformulated as follows. {\it With probability tending to 1 at } 
$ n \to \infty $
$$
|| \hat{f} - f||^2 \le B(n)/(1 - \gamma),
$$
and ACI is constructed by replacing the values $ B(n), \gamma $ by their 
consistent estimates:
$$
|| \hat{f} - f ||^2 \le \tau^*(n) \ \frac{ \tau(2M) - 2 \tau(M)}{ 3 \tau(2M) - 2 
\tau(M)  - \tau(4M) }.
$$
A more exact result will be obtained by taking into account the following term 
of the expansion of the value $ || \hat{f} - f ||^2: $
$$
\frac{ || \hat{f} - f ||^2}{ B(n) } \le \frac{1}{ 1 - \gamma} + \frac{\zeta}
{\sqrt{N^0(n)}} (1 + \epsilon_n),
$$
where $ \epsilon_n \to 0; \ {\bf P}(|\zeta| > u) \le 2 \exp( - C u^{r/2}) $ and 
$ C $ no longer depends on $ n $. Equating the probability 
$ {\bf P}(|\zeta| > u) $, more exactly its estimate 
$ 2 \exp( - C u^{r/2} ) $ to the value $ \delta, \ \delta \approx 0+, $ we will 
easily find $ u = u(\delta) $ and construct an approximate ACI with reliability 
$ \approx 1 - \delta $ of the form
$$
|| \hat{f} - f ||^2 \le \frac{ \tau^*(n)}{1 - \hat{ \gamma}} + \tau^*(n) \frac{
u(\delta) } { \sqrt{N(n)} }.
$$
Closer consideration reveals an effect that somewhat reduces the exactness of 
ACI. Let (as is true in all the three considered problems under the formulated 
assumptions)
$$
{\bf P} \left( ||\hat{f} - f||^2 / B(n) > u \right) \le \exp( 
- \phi(C_1 u)),
$$
$$
{\bf P}\left( \tau^*(n)/B(n) < 1/u \right) \le 
\exp( - \phi(C_2 u)), \ u > C,
$$
where at $ u \to \infty \ \Rightarrow \ \phi(u) \to 0. $ We denote
$$
{\bf Q}(u) = {\bf P} \left( ||\hat{f} - f||^2 / \tau^*(n) > u \right).
$$
{\bf Theorem $ \tau $.}  {\it At $ u \le C/B(n) $ the following inequality 
holds:}
$$
Q(u) \le 2 \exp( - \phi(C \sqrt{u})).
$$
{\bf Proof.} We have by the full probability formula we (we shall understood 
$ {\bf P}(A \big/ B) $ 
as  the conditional probabilities, if, of course, $ A $ and $ B $ are events):
$$
{\bf Q}(u) \le {\bf P} \left( \frac{ ||\hat{f} - f||^2}{\tau^*(n)} > u \big/
\frac{\tau^*(n)}{B(n)} > \frac{1}{v}  \right) \cdot
{\bf P} \left( \frac{\tau^*(n)}{B(n)} > \frac{1}{v} \right) +
$$

$$
+ {\bf P} \left( \frac{ ||\hat{f} - f||^2}{\tau^*(n)}> u \big/ \frac{\tau^*(n)}{B(n)} \le 
\frac{1}{v} \right) \cdot {\bf P} \left( \frac{\tau^*(n)}{B(n)} \le \frac{1}{v} \right)
\stackrel{def}{=} Q_1 + Q_2;
$$

$$
Q_1 \le {\bf P} \left( || \hat{f} - f||^2/B(n) > u/v \right) 
\le \exp( - \phi(C_1 u/v));
$$

$$
Q_2 \le {\bf P} \left( \tau^*(n)/B(n) \le 1/v  \right) \le 
\exp \left( - \phi(C_2 v) \right).
$$
 Summing up and put $  v = C_3 \sqrt{u} $, we obtain  the assertion of 
the theorem. \par
 The increase in the probability $ {\bf Q} $ compared to $ {\bf P}_f $ is 
apparently explained by 
the ability of the denominator, i.e. $ \tau^*(n) $ to take values close to 
zero. \par
 Note in conclusion that the estimates proposed by us have successfully passed 
experimental tests on problems {\bf R, D, S } by simulate  modeled with the 
use of pseudo-random numbers as well as on real data (of seismic signals etc.) 
for which our estimations of the spectrum were compared with classical estimates 
obtained by the spectral window method.\par

\vspace{2mm}

{\bf Aknowledgement.} It is the author's pleasure to convey my gratitude to I.A.Ibragimov 
and V.N.Sudakov for fruitful discussions of the problems of the problems under consideration, 
to D.Donoho for sending me his publications and manuscripts, as well as A.Pridor for the 
possibility to practically implement of our  estimations. \\  

\newpage

{\bf References} \\

\vspace{3mm}

 AAD W. VAN DER VAART and MAARK J. van Der Laan. Smooth Estimation of a monotonic density.
{\it Statistics.} 2003, V, 37, $ N^o $ 3, p. 189 - 203.\\ 

 Allal J., and Kaaouachi. Adaptive R - estimation in a Linear Regression Model with Arma 
 Errors. {\it Statistics.} 2003,  V. 37, $N^o $ 4, July - August, pp. 271 - 286.\\

 Bernstein S. {\it Collected Works.} Moscow, 1952, AN SSSR, v.1 (in Russian). \\

 Bernstein S.{\it Collected Works.} Moscow, 1954, AN SSSR, v.2 (in Russian).\\

 Bobrov P.B., Ostrovsky E.I.  Confidence intervales by adaptive estimations.
{\it  Zapiski Nauchn. seminarov POMI.}  St. - Petersburg, 1997, v. 37 b.2, 28 - 45.\\ 

 Candes E.J.,  Ridgelets: Estimating with ridge Functions.  {\it Annales of 
 Statistics,} 2003, v. 31 $ N^o $ 31, 1561 - 1569.\\

 Corrine Berzin, Jose' R. Leon and Joaquim Ortega. Convergence of non - linear 
functionals of Smoothed Empirical Processes and Kernel Density Estimates. 
{\it Statistics, } 2003, V. 37 $ N^o $ 4, pp. 217 - 242. \\

 Dette H. and Melas V. Ch. B.  Optimal Design for Estimating individual 
 coefficients in Fourier Regression Model.{\it Annales of Statistics,} 2003, v. 31
 $ N^o $ 5, 1669 - 1692.\\

 Donoho D., Jonstone I., Keryacharian G., Picard D.  Density estimation 
 by wavelet thresholding.{\it  Technical report} $ N^0 426,$  1993, Dept. of Stat., 
 Stanford University.\\ 

 Donoho D., Jonstone I.  Adapting to unknown smoothness via wavelet 
 shrinkage.  {\it Technical report} $ N^0 425,$ 1993,  Dept. of Stat.,
 Stanford University.\\

 Donoho D.  Wedgelets: nearly minimax estimation of edges.{\it  Annales of 
 of Statist.,} 1999, v. 27 b. 3 pp. 859 - 897.\\

 Donoho D.  Unconditional bases are optimal bases for data compression 
 and for statistical estimation. {\it Applied Comput. Harmon. Anal., } 1996,
 v. 3 pp. 100 - 115.\\

 Donoho D.  Unconditional bases and bit - level compression.{\it Appl. 
 Comput. Anal.,}  1999, v. 3 pp. 388 - 392.\\

 Efroimovich S.  Nonparametric estimation of the density of a unknown 
 smoothness. {\it Theory Probab. Appl.,} 1985, v. 30 b. 3,  557 - 568.\\

 Fiegel T., Hitczenko P., Jonson W.B., Shechtman G., Zinn J.  Extremal .
 properties of Rademacher functions with applications to the Khinchine and Rosental 
 inequalities. {\it Transactions of the American Math. Soc.,} 1997, v. 349 $ N^o $ 3, 
 997 - 1024.\\ 

 Golubev G., Nussbaum M.  Adaptive spline Estimations in the nonparametric 
 regression Model. {\it Theory Probab. Appl.,} 1992, v. 37 $ N^o $ 4, 521 - 529.\\

 Golybev G. Nonparametric estimation of smooth spectral densities of 
 Gaussian stationary sequences. Theory Probab. Appl., 1994, v. 38 b. 2, 28 - 45.\\

 Ibragimov I.A., Khasminsky R.Z.  On the quality boundaries of nonparametric 
 estimation of regression. {\it  Theory Probab. Appl.,} 1982, v. 21 b. 1, 81 - 94.\\

 Jonson W.B., Schechtman G., Zinn J. Best Constants in the moment Inequalities for linear
 combinations of independent and exchangeable random Variables. {\it  Annales Probab.,}
 1985, v. 13, 234 - 253.\\

 Korolyuk V.S., Borovskich Yu. P. {\it Theory of $ U - $ Statistics,} 1993, Springer,
 Berlin - Heidelberg - New York - Tokyo. \\

 Kozachenko Yu.V., Ostrovsky E.I.  Banach Spaces of random Variables of subgaussian type.A
 {\it Theory Veroyatn. Mathem. Statist.,} 1983, v. 32, 52 - 53.\\

 Lee Geunghee.  Choose of smoothing Parameters in Wavelet Series Estimators.
Journal of Nonparametric Statistics, 2003, v. 15 (4 - 5), p. 421 - 435.\\

 Lepsky O.  On adaptive estimation problem in the Gaussian white Noise.{\it Theory Probab.
 Appl.,}  1990, v. 35 b. 3, 454 - 461.\\

 Nikolsky S. {\it Inequalities for integer Functions of finite Power add Their 
 Applications in the Theory of differentiable Functions of many Variables.} (in Russian).
 In: Trudy Mathemat. Inst. im. V.V.Steklova AN SSSR, 1951, v. 51, 244 - 278.\\

 Nussbaum M.  Spline smoothing in regression Models and Asymptotic Efficiency in  $ L_2. $ 
{\it Annales of Statist.,} 1985, v. 13 b. 3, 984 - 997.\\

 Ostrovsky E.I. Adaptive estimation in three classical 
 problem  on nonparametric statistics. {\it Aktualnye problemy 
 sovremennoy mathematiki.}  Novosibirsk, NII MI OO, 1997, v.3 pp. 142 - 146.\\

 Ostrovsky E.I. The adaptive estimation in multidimensional statistics.A
 In: {\it Proseeding of the 5th international conference on simulation of devices 
 and technologies (ISDT).}  Obninsk, 1996, pp. 115 - 118.\\

 Ostrovsky E.I. {\it Exponential estimates for random fields and their 
 applications } (in Russian). Obninsk, OIATE, 1999, 350 pp.\\

 Pizier G. Condition  d'entropie  assurant la continuite de certains 
 processes et applications a  l'analyse  harmonique. In: {\it Sem. 
  d' anal. funct.,} 1979 - 1980, v. 23 - 24, 1 - 43. \\ 

 Plicusas A. Some Properties of Multiply Integral Ito.{\it Liet. Mathem. Rink.,}
 1981, v. 21 b. 2, 163 - 173.\\

 Polyak B., Tsybakov A. $ C_p - $ criterion in projective Estimation of 
 Regression.{\it  Theory Probab. Appl.,} 1990, v. 35 b. 2, 293 - 306.\\

 Polyak B., Tsybakov A. A family of asymptotically optimal Methods for selecting the 
 Order of projectiv Estimation of Regression.{\it Theory Probab. Appl.,} 1992, v. 37, b. 3, 
 471 - 485.\\

  Ronzin A.  Asymptotic formulae for moments of $ U \ - $ statistics with degenerate 
  kernel.{\it Theory Probab. Appl.,} 1982, v. 27 b. 2, 163 - 173.\\

  Rozental H. On the subspaces of $ L_p (p > 2) $  spanned by sequence of
  independent variables. {\it Probab. Theory Appl.,} 1982, v. 27 b.1, 47 - 55.\\
 
  Saulis L., Statuliavichius V. {\it Limit Theorems for Great Deviations.} Vilnius,
  Mokslas, 1989. \\

  Shiryaev A.N. {\it Probability.} Kluvner Verlag, 1986.\\ 

  Tchentsov N.N. {\it Statistical resolving rules and optimal conclusions.}
  Moscow, Nauka, 1972. 529 pp. \\

  Timan A. {\it Theory of Approximation of Functions of Real Variables } (in Russian).
  Moscow, GIFML, 1960.\\
   
  Tony Cai T.  Adaptive wavelet Estimation: a Block Thresholding and Oracle 
  Inequality Approach.{\it  Annales of Math. Statist.,} 1999, v. 27 b. 3, 898 - 924.\\

\end{document}